\documentclass[11pt]{amsart}
\usepackage[T1]{fontenc}
\usepackage{lmodern,microtype,amssymb,mathtools}
\usepackage[a4paper,hmargin=28mm,top=24mm,bottom=25mm]{geometry}
\usepackage{xcolor}
\usepackage[colorlinks=true,linkcolor=blue!45!black,citecolor=blue!45!black,urlcolor=blue!45!black]{hyperref}
\hypersetup{pdftitle={Topological obstructions to embedding a matrix algebra bundle into a trivial one},pdfauthor={A. V. Ershov}}
\newtheorem{theorem}{Theorem}[section]
\newtheorem{prop}[theorem]{Proposition}
\newtheorem{lemma}[theorem]{Lemma}
\newtheorem{cor}[theorem]{Corollary}
\theoremstyle{definition}\newtheorem{example}[theorem]{Example}
\theoremstyle{remark}
\newcommand{\C}{\mathbb C}\newcommand{\Z}{\mathbb Z}\newcommand{\Q}{\mathbb Q}
\newcommand{\one}{\mathbf 1}\newcommand{\triv}[1]{\underline{\mathbb C}^{\,#1}}
\newcommand{\alg}[1]{\underline{M}_{#1}}\newcommand{\PU}{\mathrm{PU}}\newcommand{\SU}{\mathrm{SU}}
\DeclareMathOperator{\End}{End}\DeclareMathOperator{\Aut}{Aut}\DeclareMathOperator{\Emb}{Emb}
\DeclareMathOperator{\Gr}{Gr}\DeclareMathOperator{\rk}{rk}
\DeclareMathOperator{\Br}{Br}\DeclareMathOperator{\im}{im}\DeclareMathOperator{\Ind}{Ind}
\newcommand{\dd}{\delta}\newcommand{\id}{\mathrm{id}}\newcommand{\red}{\rho_k}
\newcommand{\FinBr}{\mathcal B_{\mathrm{emb}}}
\title[Embeddings of matrix algebra bundles]{Topological obstructions to embedding a matrix algebra bundle into a trivial one}
\author{A. V. Ershov}
\date{Revised manuscript, 22 November 2026}
\subjclass[2020]{Primary 55R40, 19L20; Secondary 55S35, 46L80}
\keywords{Matrix algebra bundle, unital embedding, obstruction theory, topological Brauer group, UHF algebra}
\begin{document}
\begin{abstract}
Let $A$ be a bundle of complex $k\times k$ matrix algebras over a connected finite CW complex $X$. We characterize the existence of a unital embedding into a trivial $M_{kt}$-bundle for some $t$ coprime to $k$. Its Brauer class must vanish, and for $A=\End(V)$ the necessary and sufficient condition is $[V]\in kK^0(X)$. If it holds, every multiplicity $t$ with $2t+1\ge\dim X$ works for $k\ge2$. Through dimension five we give an equivalent criterion in characteristic classes for coprime multiplicities $t\ge2$ and compute the indeterminacy of the second obstruction. On $S^{2r}$ with $k\ge r\ge2$ we determine all admissible multiplicities, including those below the stable range for the multiplicity bundle. We also identify the kernel of the map from the group obtained by quotienting matrix bundles by coprime-embeddable factors to the ordinary Brauer group. The resulting exact sequence has no natural splitting, and a Moore-space example rules out representability of this quotient on finite complexes. Finite frames give integral clutching units for $KU$-lines on suspensions. Along a fixed coprime tensor-power tower, the stabilized frames act up to coherent homotopy on the stabilized matrix Grassmannian.
\end{abstract}
\maketitle
\section{The embedding problem}
Given a locally trivial bundle $A\to X$ with fiber $M_k=M_k(\C)$, for which positive integers $t$ is there a unital embedding
\begin{equation}\label{eq:embedding}
                      A\longrightarrow \alg{kt}=X\times M_{kt}\,?
\end{equation}
We call $t$ the \emph{multiplicity}. All algebra maps preserve the involution, and embeddings and trivializations are continuous over the indicated base. We assume that $X$ is a connected finite CW complex, put $d=\dim X$, and take $k\ge2$ unless stated otherwise. The case $k=1$ has the scalar embedding for every $t$. Complex vector bundles carry Hermitian metrics. An isomorphism of such bundles can be made unitary by polar decomposition.

Even over $S^2$, the answer depends on the arithmetic of the target. If $H$ is the Hopf line, set $V=H\oplus\triv{k-1}$. The algebra $\End(V)$ has zero Brauer class. An inclusion into $\alg{kt}$ would yield a rank-$t$ bundle $W$ with $V\otimes W$ trivial. Taking determinants gives $t c_1(H)+k c_1(W)=0$, so $k\mid t$. Conversely, if $k\mid t$ one can choose $c_1(W)=-(t/k)c_1(H)$; over $S^2$ rank and first Chern class classify complex bundles. Thus even allowing all coprime multiplicities leaves an obstruction. Over higher spheres there is a further restriction coming from the geometric rank of $W$.

Write $R=K^0(X)$, $I=\ker(\rk:R\to\Z)$, and regard integers as multiples of the unit of $R$. Let $\dd(A)\in H^3(X;\Z)$ denote the obstruction to lifting the projective unitary frame bundle of $A$ to a unitary bundle. Our main criterion is the following.
\begin{theorem}\label{thm:main}
For an $M_k$-bundle $A$ the following are equivalent.
\begin{enumerate}
\item There is an embedding \eqref{eq:embedding} at some multiplicity coprime to $k$.
\item $\dd(A)=0$ and, writing $A\cong\End(V)$, the class $[V]$ is divisible by $k$ in integral $K^0(X)$.
\item There is an embedding \eqref{eq:embedding} for every positive $t$ with $2t+1\ge d$.
\end{enumerate}
The divisibility condition has the same truth value for every choice of $V$. For a specified smaller multiplicity the exact test, in the split case, is the existence of $W$ with $V\otimes W\cong\triv{kt}$.
\end{theorem}

The question and its matrix Grassmannian model were studied in \cite{Ershov2000,Ershov2008,v14,published}. The stable relation between embedded bundles and rank-one virtual bundles already appears in \cite{Ershov2000}. Here the target size is kept explicit. The proof of Theorem~\ref{thm:main} uses multiplication in $R/kR$; it imposes no lower bound on the degree $k$ of the given algebra. Theorem~\ref{thm:four} identifies the degree-four obstruction as a coset and makes the dependence on a special-unitary lift explicit. Theorem~\ref{thm:spheres} gives the complete answer on $S^{2r}$ when $k\ge r\ge2$, with no restriction on the target multiplicity. Sections~\ref{sec:uhf} and~\ref{sec:quotient} treat, respectively, passage up a specified UHF tower and the group defined by tensoring with coprime-embeddable factors. Proposition~\ref{prop:lowbr} identifies this quotient with $H^2(X;\Q/\Z)$ through dimension three. Proposition~\ref{prop:moore} uses two Moore spaces to test natural splitting and representability. Section~\ref{sec:fred} derives integral suspension twists from the family index and identifies the homotopy quotient of the stabilized Grassmannian by the stabilized frames.

\section{Frames, centralizers and multiplicity bundles}\label{sec:frames}
Give $F_{k,t}=\Emb(M_k,M_{kt})$ the point-norm topology. Its quotient by precomposition with $P=\PU(k)$ is the matrix Grassmannian of images. Unitary conjugacy of systems of matrix units gives
\[
 F_{k,t}=\PU(kt)/(1_k\otimes\PU(t)),\qquad
 \Gr_{k,t}=\PU(kt)/(\PU(k)\otimes\PU(t)).
\]
Thus $F_{k,t}\to\Gr_{k,t}$ is a principal $P$-bundle. For $(k,t)=1$ there is also an identification $\Gr_{k,t}=\SU(kt)/(\SU(k)\times\SU(t))$ with the tensor-product inclusion. Indeed its kernel is trivial and the centers of the two factors generate the whole center of $\SU(kt)$. In particular, the tautological $M_k$-bundle on $\Gr_{k,t}$ has an $\SU(k)$-lift, given by the principal bundle with total space $\SU(kt)/(1_k\otimes\SU(t))$.

\begin{prop}\label{prop:grass}
Let $P_A\to X$ be the bundle of algebra frames of $A$. Embeddings \eqref{eq:embedding} are sections of $P_A\times_P F_{k,t}$. The associated universal bundle over $BP$ has total space homotopy equivalent to $\Gr_{k,t}$. Its projection corresponds to the classifying map for $F_{k,t}\to\Gr_{k,t}$.
\end{prop}
\begin{proof}
In a local frame a section is a family of embeddings of $M_k$; the change of local frame is precomposition. The map $EP\times_P F_{k,t}\to\Gr_{k,t}$ is the associated bundle with fiber $EP$, which is contractible. Pulling back $F_{k,t}$ and $EP$ to the total space gives isomorphic principal bundles, with common model $EP\times F_{k,t}$. This identifies the projections up to homotopy and gives the lifting criterion. Over $BU(k)$ the corresponding universal total space is the homotopy pullback along $BU(k)\to BP$. For coprime $k,t$ it is homotopy equivalent to $\Gr_{k,t}\times K(\Z,2)$: the pulled-back central lifting fibration is principal with fiber $K(\Z,2)$ and the tautological lift gives it a section.
\end{proof}

A fiberwise unital embedding has a commutant $C$ of degree $t$, and multiplication identifies $A\otimes C$ with $\alg{kt}$. These statements follow by choosing matrix units locally. Local changes of an embedded frame can be implemented by continuous unitaries, either from the homogeneous-space description or by the polar decomposition of an intertwiner of close frames. Recall the identities $k\dd(A)=0$ and $\dd(A\otimes C)=\dd(A)+\dd(C)$. They also follow by multiplying the scalar triple-overlap cocycles of local unitary lifts and taking determinants. Hence the Brauer class of an embedded $A$ is killed by both $k$ and $t$. If $(k,t)=1$ then $A\cong\End(V)$ for a rank-$k$ bundle $V$. Two such bundles differ by tensoring with a line: the module homomorphisms between their two defining modules form a line bundle, and evaluation gives the isomorphism.

\begin{lemma}[Multiplicity bundle]\label{lem:mult}
For $A=\End(V)$, giving an embedding \eqref{eq:embedding} is equivalent to giving a Hermitian rank-$t$ bundle $W$ and a unitary isomorphism $V\otimes W\cong\triv{kt}$, up to a compatible unitary isomorphism of $W$. Its commutant identifies with $\End(W)$.
\end{lemma}
\begin{proof}
An embedding makes $E=\triv{kt}$ a left $\End(V)$-module. Put $W=\operatorname{Hom}_{\End(V)}(V,E)$. In a local frame of $V$, evaluation on its first basis vector identifies this space with the range of the projection associated to $e_{11}$. It is therefore a bundle of constant rank $t$. Evaluation $V\otimes W\to E$ is an isomorphism in every fiber. The invariant metric on the multiplicity space makes evaluation unitary. Conversely transport the action $a\mapsto a\otimes 1_W$ through the chosen trivialization. The two constructions are inverse in the stated sense; the commuting endomorphisms act on the multiplicity factor.
\end{proof}

Over a field there is a useful arithmetic comparison. Suppose that a central simple algebra $A$ of degree $k$ over a field $F$ embeds into $M_{kt}(F)$ with $(k,t)=1$. Write $A=M_s(\Delta)$, where $\Delta$ is a central division algebra of degree $e$, so $k=se$. Every finite left $A$-module is a direct sum of copies of the simple module $\Delta^s$, of $F$-dimension $se^2=ke$. Applying this to $F^{kt}$ shows that $e$ divides $t$. Since $e$ also divides $k$, $e=1$. For bundles an analogous argument kills the Brauer class, but the module trivialization in Lemma~\ref{lem:mult} still carries topological information.

\subsection*{Integral divisibility and the prescribed target}
We use two standard facts about complex $K$-theory \cite{Husemoller,HatcherVB}. The rank-zero ideal $I$ is nilpotent, by the multiplicative skeletal filtration; in particular every $1+z$, $z\in I$, is a unit. A virtual class of positive rank $s$ has a representative of that rank if $d\le2s+1$, and bundles of rank $s$ with the same stable class are isomorphic if $d\le2s$. The latter two assertions follow from the $(2s+1)$-connectivity of $BU(s)\to BU$, using $U(s+1)/U(s)=S^{2s+1}$. The one-dimension difference matters only in using the existence statement for $W$.

\begin{proof}[Proof of Theorem~\ref{thm:main}]
An embedding at a coprime multiplicity gives $A=\End(V)$ by the commutant argument. By Lemma~\ref{lem:mult},
\[
                         [V][W]=kt\quad\text{in }R.
\]
In $R/kR$ the class of $W$ is a unit: its rank $t$ is a unit modulo $k$ and its reduced part is nilpotent. Thus $[V]=0$ in $R/kR$. This is precisely integral divisibility, even when $R$ has torsion. Replacing $V$ by $V\otimes L$ multiplies its class by the unit $[L]$, proving independence of the lift.

Conversely, write $[V]=ku$. The element $u$ has rank one, hence has an inverse in $R$. Choose a rank-$t$ representative $W$ of $t u^{-1}$ when $d\le2t+1$. The class of $V\otimes W$ is $kt$. Since $k\ge2$, $d\le2t+1\le2kt$, so cancellation gives an actual trivialization. Lemma~\ref{lem:mult} constructs the embedding. Arbitrarily large integers coprime to $k$ exist, proving the remaining implication.
\end{proof}

\begin{cor}\label{cor:generalrank}
Suppose $A=\End(V)$ and $d\le2t+1$. For a prescribed $t$, with arbitrary $(k,t)$, an embedding exists exactly when $kt\in [V]R$.
\end{cor}
\begin{proof}
A solution $[V]w=kt$ has rank $\rk(w)=t$. The same realization and cancellation argument applies. Necessity holds in every dimension.
\end{proof}

\section{The first two obstructions}\label{sec:four}
Write $s(A)\in H^2(X;\Z/k)$ for the lifting class of the central extension $\mu_k\to\SU(k)\to\PU(k)$. Fix the choice of generator for which $s(\End(V))=\red c_1(V)$, where $\red$ denotes reduction modulo $k$. Indeed twisting $V$ by a line changes $c_1(V)$ by $k$ times its first Chern class, and the determinant induces this reduction on the universal class. The integral Bockstein of $s(A)$ gives $\dd(A)$ with compatible conventions. An embedding at coprime multiplicity implies $s(A)=0$ by the tautological special-unitary lift of Section~\ref{sec:frames}.

Suppose $s(A)=0$, and choose an $\SU(k)$-lift with associated vector bundle $V$. For $T_k(X)=H^2(X;\Z)[k]$ define
\begin{equation}\label{eq:q}
 q_k:T_k(X)\longrightarrow H^4(X;\Z/k),\qquad
 q_k(a)=\red\!\left(\binom{k}{2}a^2\right).
\end{equation}
This map is additive because the coefficient of $ab$ in its polarization is $k(k-1)$. Its image is killed by $2$ and is zero when $k$ is odd.
\begin{theorem}\label{thm:four}
When $s(A)=0$, the class
\[
 o_4(A)=\big[\red c_2(V)\big]\ \in\ H^4(X;\Z/k)/q_k(T_k(X))
\]
is independent of the chosen lift and natural under pullback. For $k,t\ge2$ coprime and $d\le5$, the two successive conditions $s(A)=0$, $o_4(A)=0$ are necessary and sufficient for an embedding into $\alg{kt}$. They are necessary in every dimension.
\end{theorem}
\begin{proof}
With the identification of the projective bundle fixed, lifts through the central covering form an $H^1(X;\Z/k)$-torsor. On vector bundles the change is $V\mapsto V\otimes L$, with a specified trivialization of $L^{\otimes k}$. The Bockstein sequence shows that the possible $a=c_1(L)$ exhaust $T_k(X)$. The splitting principle gives
\[
 c_2(V\otimes L)=c_2(V)+(k-1)a c_1(V)+\binom{k}{2}a^2.
\]
Since $c_1(V)=0$, the possible residues form exactly the coset in the statement. This construction commutes with pullback.

For necessity, assume there is an embedding. Choose a special-unitary lift $V$ and its multiplicity bundle $W$. Taking the determinant of the trivial $V\otimes W$ gives $k c_1(W)=0$. Multiplication by $t$ is an automorphism of $T_k(X)$. Choose $a\in T_k(X)$ with $ta=c_1(W)$ and a line $L$ of class $a$. Replace $(V,W)$ by $(V\otimes L,W\otimes L^{-1})$. Both determinants are now trivializable. For this pair the degree-four tensor-product formula gives
\begin{equation}\label{eq:c2}
                         0=t c_2(V)+k c_2(W),
\end{equation}
where the notation refers to the adjusted bundles. This identity holds integrally, including on torsion, by the splitting principle. Coprimality gives $\red c_2(V)=0$ for this lift.

Conversely, vanishing of the coset supplies an $\SU(k)$-lift with $c_2(V)=k b$ for some $b\in H^4(X;\Z)$. The map $c_2:B\SU(n)\to K(\Z,4)$ is 5-connected for $n\ge2$: the source is 3-connected, $c_2$ induces an isomorphism on $\pi_4$, and $\pi_5 K(\Z,4)=0$. It is 6-connected for $n\ge3$, since also $\pi_5 B\SU(n)=\pi_4\SU(n)=0$ by Bott periodicity in the stable range. Thus through dimension five every fourth cohomology class is realized at rank $t\ge2$, and $c_2$ classifies special-unitary bundles of rank at least three. Choose an $\SU(t)$-bundle $W$ with $c_2(W)=-tb$. The tensor product has trivial determinant and zero second Chern class by the tensor-product formula used in \eqref{eq:c2}. Its rank is $kt\ge4$, so it is trivial. Apply Lemma~\ref{lem:mult}.
\end{proof}

The dimension bound is sharp. On $S^6$ take $k=5$, $t=3$ and $[V]=5+b$, where $b$ is a Bott generator. The groups containing both obstructions vanish. A rank-three multiplicity bundle with class $3+nb$ would give $5n+3=0$ in $\Z$, which is impossible since $b^2=0$.

\begin{example}[Sharpness of the indeterminacy]\label{ex:lens}
For every even $k$, let $X$ be the four-skeleton of the infinite lens space $B\mu_k$, and let $L$ come from its defining character. For $a=c_1(L)$ the classes $a$ and $a^2$ generate $H^2(X;\Z)=\Z/k$ and $H^4(X;\Z)=\Z/k$ respectively. Indeed the circle fibration $S^\infty/\mu_k\to\mathbb{CP}^{\infty}$ has Euler class $k$ times a generator, giving the ring $\Z[a]/(ka)$ by Gysin. In the standard cellular model the odd cellular boundaries are zero, so truncation at degree four preserves these two groups \cite{HatcherAT}. The two lifts $\triv{k}$ and $L^{\oplus k}$ of the trivial algebra have second Chern residues differing by the nonzero element of order two, $\binom{k}{2}\red(a^2)$. Thus even the zero algebra class can have a nonzero residue for a particular lift.
\end{example}

For a closed-manifold example in the strict range $d<2\min(k,t)$, take $(X,k,t)=(\mathbb{RP}^4,6,5)$ and let $L$ be the complexified tautological line. The mod-two reduction of $c_1(L)$ is $u^2$, where $u$ generates $H^1(X;\Z/2)$, because the underlying real bundle is two copies of the tautological real line. Hence $15c_1(L)^2$ is already nonzero modulo two. The lifts $L^{\oplus6}$ and $\triv{6}$ have the same endomorphism algebra and give different residues modulo six. This is the ambiguity in the statement preceding Theorem~8 of \cite{v14}.

\begin{example}[A class killed by amplification]\label{ex:coefficient}
Take $X=\mathbb{RP}^2$ and $V=\one\oplus L$ for the same complexified line. The first obstruction of $\End(V)$ is nonzero, whereas $c_1(V^{\oplus2})=0$. Thus tensoring the algebra with $M_2$ kills the primary obstruction. Indeed $H^2(X;C)=C/2C$ for any coefficient group $C$, and the coefficient inclusion $\Z/2\to\Z/4$, $1\mapsto2$, induces the zero map in degree two.
\end{example}

\section{Higher obstructions in the stable range}\label{sec:higher}
This section identifies the cohomological stages of the original lifting problem. They give successive tests when information on the skeleta of $X$ is available. Fix $k,t\ge2$ coprime. In the range $1\le j\le2t$, the homotopy sequence of
\[
            \PU(t)\xrightarrow{1_k\otimes -}\PU(kt)\longrightarrow F_{k,t}
\]
and Bott periodicity give $\pi_jF_{k,t}=\Z/k$ for odd $j$ and zero for even $j$. On stable odd homotopy in degrees at least three the first arrow is multiplication by $k$. On fundamental groups it is the injection $\Z/t\to\Z/(kt)$, $1\mapsto k$. At the endpoint, $\pi_{2t}\PU(kt)=0$ and the map on $\pi_{2t-1}$ from $\PU(t)$ is injective, so $\pi_{2t}F_{k,t}=0$ without requiring stable homotopy for $\PU(t)$ in degree $2t$. The orbit map $\PU(k)\to F_{k,t}$ factors through $u\mapsto u\otimes1_t$ in $\PU(kt)$. On $\pi_1$, and in odd degrees $2r-1$ with $2\le r\le\min(k,t)$, it sends a generator to $t\bmod k$, up to sign. We normalize the coefficient generator in these degrees by this unit.

\begin{prop}\label{prop:tower}
If $d\le2t+1$, an embedding exists exactly when the relative Moore--Postnikov tower of $\Gr_{k,t}\to B\PU(k)$ admits compatible partial lifts of the classifying map of $A$ through the stages required in dimensions $2r\le d$. On a specified preceding lift the obstruction at that stage belongs to $H^{2r}(X;\Z/k)$. The first stage is $s(A)$; on a specified $\SU(k)$-lift the next is $\red c_2(V)$ in the coefficient convention above.
\end{prop}
\begin{proof}
Replace the Grassmannian map by the fibration in Proposition~\ref{prop:grass}. A nonzero $\pi_{2r-1}$ of its fiber gives a relative Postnikov stage with obstruction in degree $2r$. The base is simply connected. The first stage can be identified over $B\PU(k)$ with $B\SU(k)$: the tautological special-unitary lift gives the comparison, and the boundary on $\pi_2$ identifies the first fiber group with $\mu_k$. This stage is simply connected, and all later stage fibers have degree at least three, so the subsequent coefficient systems are constant. A vanishing obstruction permits a choice of lift, equivalently a null-homotopy of the pulled-back class. When this is possible, the choices up to fiberwise homotopy over a fixed preceding map form a torsor for $H^{2r-1}(X;\Z/k)$. Cellular lifting terminates at $d$. This is the relative Moore--Postnikov construction \cite{Whitehead}.

The chosen tautological lift factors the Grassmannian map as $\Gr_{k,t}\to B\SU(k)\to B\PU(k)$. The homotopy fiber of the first arrow is the total space $\SU(kt)/(1_k\otimes\SU(t))$ of the tautological $\SU(k)$-bundle. It is the simply connected cover of $F_{k,t}$. By comparison, base change of the original fibration along $B\SU(k)\to B\PU(k)$ still has fiber $F_{k,t}$. The action of $\SU(k)$ by changing the frame has orbit map given by $u\otimes1_t$ followed by projection. On $\pi_3$ this is multiplication by $t$ followed by reduction modulo $k$. Evaluation on the clutching sphere therefore identifies the degree-four class with $\red c_2$, with the stated normalization. This is the chosen-lift version of Theorem~\ref{thm:four}.
\end{proof}

\begin{prop}[Divided Chern classes]\label{prop:divided}
Let $2\le r\le\min(k,t)$ and let $Y_r=BU(k)\langle2r\rangle$ be the $(2r-1)$-connected cover with projection $i_r$ to $BU(k)$. For the universal embedding problem pulled back to $Y_r$, the first possible obstruction is the reduction of
\[
                  \widetilde c_r=i_r^*c_r/(r-1)!\in H^{2r}(Y_r;\Z),
\]
with compatible choices of Bott and obstruction generators. The same formula holds on pulling back along a specified map $X\to Y_r$ with its induced lower-stage data.
\end{prop}
\begin{proof}
In this range $2r-1<2k$, so $\pi_{2r}(Y_r)=\pi_{2r-1}U(k)=\Z$. By Hurewicz and universal coefficients, an integral class in degree $2r$ is determined by its evaluation on a generator. Bott integrality gives the value $\pm(r-1)!$ for $i_r^*c_r$, proving both divisibility and primitivity of the divided class \cite[Ch.~20, Cor.~9.8]{Husemoller}. All lower reduced cohomology groups on $Y_r$ vanish, so the lower obstructions and their lift indeterminacies vanish there. Evaluate on a generating $S^{2r}$. The clutching map into $U(k)$ represents a stable generator; its orbit in $F_{k,t}$ represents $\pm t\bmod k$ by the computation above. This determines the universal obstruction in $H^{2r}(Y_r;\Z/k)$ and proves the claimed normalization. Naturality gives the last assertion.
\end{proof}

For instance, on $S^8$ take $k=6$, $t=5$, and let $[V]=6+b$ for a Bott generator $b$. Then $c_4(V)$ is zero modulo six since its integral value is $\pm 6$ times a generator. The divided class detects the obstruction. The easier test $[V]\notin6K^0(S^8)$ gives the same answer by Theorem~\ref{thm:main}.

The restriction on dimensions keeps the calculation at a finite target. For instance $\pi_5(\PU(2))=\pi_5(S^3)=\Z/2$ while $\pi_5(\PU(10))=\Z$. The long exact sequence gives an injection $\Z\hookrightarrow\pi_5(F_{5,2})$. Outside the range of Proposition~\ref{prop:tower} one may still use the actual relative Postnikov tower, with its actual fiber groups. The integral criterion in Theorem~\ref{thm:main} remains valid for eventual coprime embedding in every finite dimension.

\section{All target sizes on even spheres}\label{sec:spheres}
Let $X=S^{2r}$ with $r\ge2$, and choose a Bott generator $b$. Then $K^0(X)=\Z\oplus\Z b$ and $b^2=0$. Since $H^3(X;\Z)=0$ every matrix algebra bundle is split. Since $H^2(X;\Z)=0$, its defining vector bundle is unique up to isomorphism. In the following theorem the source rank is in the stable range, while the target multiplicity is arbitrary.
\begin{theorem}\label{thm:spheres}
Let $k\ge r$ and let $V$ have class $k+mb$, with $m\in\Z$. If $m=0$, then $\End(V)$ embeds at every positive multiplicity. If $m\ne0$, the exact criterion at multiplicity $t$ is
\[
                          t\ge r\qquad\text{and}\qquad k\mid tm.
\]
Writing $h=k/\gcd(k,m)$, the least multiplicity is $h\lceil r/h\rceil$. In a tower $t=l^N$ with $l\ge2$ coprime to $k$, a nonzero $m$ admits a stage precisely when $k\mid m$; the first such stage is $N=\lceil\log_l r\rceil$.
\end{theorem}
\begin{proof}
If an embedding exists, the bundle in Lemma~\ref{lem:mult} has a stable class $[W]=t+nb$ and
\[
                       (k+mb)(t+nb)=kt+(kn+tm)b.
\]
Hence $kn=-tm$. When $t<r$, the top Chern class $c_r(W)$ vanishes for rank reasons. By Bott integrality $c_r(nb)=\pm(r-1)!n$ times a generator in $H^{2r}(S^{2r};\Z)$, and therefore $n=0$. This forces $m=0$. If $t\ge r$ and $k\mid tm$, a rank-$t$ bundle realizes $t-(tm/k)b$; its tensor product with $V$ is stably trivial and cancels in rank $kt$. Finally $m=0$ makes $V$ trivial by cancellation at rank $k\ge r$. These conditions give both minimum formulas.
\end{proof}
For $X=S^6$, $k=5$, $[V]=5+5b$, every $t\ge3$ works and $t=1,2$ fail. In the dyadic tower this means that the first target is $M_{20}$. Taking $k\ge r$ and $m=k$ shows in every even dimension $2r$ that the uniform bound of Theorem~\ref{thm:main} cannot be decreased from $r$ to $r-1$.

When $2\le k<r$, a rank-$k$ bundle $V$ on $S^{2r}$ has $c_r(V)=0$, and Bott integrality gives $[V]=k$. Thus every $t$ with $kt\ge r$ works: take $W=\triv{t}$ and cancel the stable trivialization of $V^{\oplus t}$. Smaller targets can depend on the unstable bundle itself. For instance, on $S^6$ the nonzero clutching element of $\pi_5 U(2)=\pi_5 S^3=\Z/2$ gives a nontrivial rank-two bundle $V$. Its algebra cannot embed at $t=1$, since such an embedding would make $V$ a sum of two copies of a line and every line on $S^6$ is trivial. It does embed at every $t\ge2$ by the preceding argument.

\section{Passage to a UHF algebra}\label{sec:uhf}
For any $l\ge2$ put $L=M_{l^\infty}=\bigotimes_{j\ge1}M_l$ and $D=M_k\otimes L$. The unital inclusions of the stages $D_N=M_{kl^N}$ are $a\mapsto a\otimes1_l$. Bundle triviality in this section uses the point-norm automorphism group of the unital fiber.
\begin{prop}\label{prop:uhf}
An $M_k$-bundle $A$ over a finite CW complex embeds in $\alg{kl^N}$ for some $N$ if and only if $A\otimes L$ is a trivial $D$-bundle. This equivalence holds without any coprimality hypothesis. If $(k,l)=1$, it is also equivalent to condition (2) of Theorem~\ref{thm:main}.
\end{prop}
\begin{proof}
If $A$ embeds at stage $N$, let $C$ be the degree-$l^N$ commutant. The bundle $C\otimes L$ has fiber isomorphic to $L$ and is trivial, since $\Aut(L)$ is contractible \cite[Thm.~2.3]{DP}. Tensoring $A\otimes C\cong\alg{kl^N}$ by $L$ now gives the asserted trivialization of $A\otimes L$.

Conversely, a trivialization restricts to an embedding $A\to X\times D$. We give the finite-stage argument since the quantifier on $N$ is essential. Write $F_N=\Emb(M_k,D_N)$, $F_\infty=\Emb(M_k,D)$. There is a $\PU(k)$-equivariant map
\begin{equation}\label{eq:telescope}
                   \operatorname{Tel}_N F_N\longrightarrow F_\infty
\end{equation}
induced by the inclusions. It is a weak equivalence. Here are details of the relative approximation. For a compact family $f$ in $F_\infty$, apply the tail trace expectation $D\to D_N$ to $f(e_{ij})$, obtaining $a_{ij}$. These converge uniformly to the given matrix units. For large $N$, let $p$ be the spectral projection of $a_{11}$ above $1/2$, let $P=\operatorname{diag}(p,\ldots,p)$, and form the row $S=(a_{i1}p)_{i=1}^k$. The original row $(f(e_{i1}))$ has initial projection $\operatorname{diag}(f(e_{11}),\ldots,f(e_{11}))$ and final projection $1$. Consequently $S^*S$ is uniformly close to $P$ and $SS^*$ to $1$. Inverting the positive part on the corner $P$ gives a row $U=S(S^*S)^{-1/2}$ with $U^*U=P$ and $UU^*=1$. Its entries give exact matrix units $g(e_{ij})=U_i U_j^*$ in $D_N$, continuously in the parameters. These formulas fix any matrix units already in $D_N$.

For $g$ sufficiently close to $f$, the polar unitary of $\sum_i g(e_{i1})f(e_{1i})$ is close to $1$ and intertwines $f$ and $g$. Its principal logarithm gives a homotopy, stationary wherever $g=f$. Applying this construction to spheres and to disks relative to their boundaries proves the weak equivalence in \eqref{eq:telescope}; compact subsets of the telescope have bounded height and can be moved into a single stage. No equivariance of the correction formulas is needed.

Apply the associated-bundle construction to \eqref{eq:telescope} using the frame bundle $P_A$. Both projections to $X$ are fiber bundles. The induced map of total spaces is a weak equivalence by the homotopy sequences. Since $X$ is a CW complex, the existing section on the right lifts up to homotopy to a map into the telescope bundle. Its projection to $X$ is homotopic to the identity, so the covering homotopy property straightens it to a section. Its compact image has bounded telescope height. Retraction of a finite subtelescope to its last stage is equivariant and gives a section with fiber $F_N$, as required. For coprime $k,l$, apply Theorem~\ref{thm:main} and take $N$ sufficiently large.
\end{proof}

\section{The quotient by coprime-embeddable factors}\label{sec:quotient}
Here we study the geometric equivalence relation of \cite[§2.4]{v14} directly on a finite complex. Call a bundle of degree $k$ \emph{coprime-embeddable} if it has an embedding \eqref{eq:embedding} at some $t$ coprime to $k$. Include the scalar algebra. On matrix algebra bundles of all degrees put $A\sim B$ if
\[
                        A\otimes C\cong B\otimes E
\]
for some coprime-embeddable $C,E$. Write $\FinBr(X)$ for the quotient. Theorem~\ref{thm:main} shows that the tensor product of two coprime-embeddable bundles is again coprime-embeddable: their defining vector-bundle classes have the form $ku,mv$, with $u,v$ of rank one. In particular $\sim$ is transitive. All trivial matrix bundles are allowed as factors.

The quotient is an abelian group. Indeed every algebra bundle admits a tensor complement to a trivial algebra. One proof is to take $l=k$ in Proposition~\ref{prop:uhf}; the fiber $M_k\otimes M_{k^\infty}$ is isomorphic to $M_{k^\infty}$ and its automorphism group is contractible. A finite-stage embedding then supplies the commutant complement. There is also a finite-dimensional proof. For a vector bundle $V$ of rank $s$, put $e=[V]-s$. As soon as $e^j=0$ the integral identity
\[
 (s+e)\sum_{i=0}^{j-1}(-e)^i s^{j-1-i}=s^j
\]
holds in $K^0(X)$, including its torsion. For $s\ge2$ increase $j$ until the sum is realized by a bundle $W$ of rank $s^{j-1}$ and stable cancellation gives $V\otimes W\cong\underline{\C}^{s^j}$. For $s=1$ use the dual line. Take $V$ to be the underlying vector bundle of $A$. Since $A\otimes A^{\mathrm{op}}\cong\End(V)$, the algebra $A^{\mathrm{op}}\otimes\End(W)$ is a tensor complement of $A$. The vector-bundle argument is due to Donovan and Karoubi \cite[Lemma~12]{DK}.

Let $R_{\Q}=R\otimes\Q$, $I_{\Q}=I\otimes\Q$, and let $U_1(R)=1+I$ and $U_1(R_{\Q})=1+I_{\Q}$ carry multiplication. Let $\Br(X)$ be the ordinary complex topological Brauer group, identified on a finite CW complex with $H^3(X;\Z)_{\rm tors}$ \cite{DK,DPBrauer}.
\begin{theorem}\label{thm:brauer}
There is a natural short exact sequence of abelian groups
\begin{equation}\label{eq:brauer}
 0\longrightarrow \frac{U_1(R_{\Q})}{\im U_1(R)}
   \xrightarrow{\ \Phi\ }\FinBr(X)\xrightarrow{\ \dd\ }\Br(X)
   \longrightarrow0 .
\end{equation}
The kernel is abstractly isomorphic to $(\Q/\Z)^b$, where $b=\dim_{\Q} I_{\Q}$. In particular, this sequence splits as a sequence of abstract abelian groups for each $X$.
\end{theorem}
\begin{proof}
All allowed factors are endomorphism bundles by Theorem~\ref{thm:main}, so the Brauer class descends. It is onto because each classical Brauer class has a finite matrix representative. The kernel consists of classes represented by $\End(V)$.

Given $v\in U_1(R_{\Q})$, clear denominators and choose a sufficiently large $n$ such that $nv$ has an integral lift represented by a rank-$n$ bundle $V$. Send $v$ to $[\End(V)]$. If $V'$ of rank $n'$ is another choice, $n'[V]-n[V']$ is torsion in $R$. Multiplication by a large common integer kills this torsion and brings the common rank into the cancellation range. Consequently suitable direct sums of $V$ and $V'$ are isomorphic. On endomorphism algebras these direct sums correspond to tensoring with constant matrices, so the construction is well-defined. Tensor products show that it is a homomorphism. Every element in the kernel of $\dd$ occurs by taking $v=[V]/\rk(V)$. Pullback gives the same construction, including when a larger rank is chosen on the new base.

If $v$ is the image of a rank-one integral class $u$, represent $nu$ itself by $V$ at sufficiently large rank. Then $\End(V)$ is coprime-embeddable. Conversely, suppose $[\End(V)]=0$ in $\FinBr(X)$. There are coprime-embeddable factors $\End(E),\End(F)$ with $\End(V\otimes E)\cong\End(F)$. Hence $V\otimes E\cong F\otimes L$ for a line $L$. By Theorem~\ref{thm:main}, $[E]=e u$ and $[F]=f w$ with $e=\rk E$, $f=\rk F$, and $u,w\in U_1(R)$. Since $f=ne$, division in $R_{\Q}$ gives $[V]/n=w[L]u^{-1}$. This proves exactness.

Finally the finite logarithm identifies $1+I_{\Q}$ with the additive rational vector space $I_{\Q}$. Under this identification the image of $1+I$ is a full lattice. To see discreteness, put $\Lambda=\im(I\to I_{\Q})$, a finitely generated lattice. The finitely many denominators in the logarithm give a fixed $M$ such that $\log(\im(1+I))\subseteq M^{-1}\Lambda$. To see full rank, for $a\in\Lambda$ choose an integer $D>0$ clearing all factorial denominators in the finite exponential of $Da$. Then $\exp(Da)\in 1+\Lambda$, so $Da$ belongs to the logarithmic image. The rational vector space modulo a full lattice is abstractly $(\Q/\Z)^b$. A divisible abelian group is injective, giving an abstract splitting of \eqref{eq:brauer}.
\end{proof}

In particular $\FinBr(S^{2r})=\Q/\Z$ for $r\ge1$, while $\FinBr(S^{2r+1})=0$ for $r\ge0$. The difference from the Dadarlat--Pennig UHF Brauer group is already visible on $S^2$. For $L=M_{l^\infty}$ one has $\pi_1\Aut(L\otimes\mathcal K)=K_1(L)=0$, where $\mathcal K$ denotes the compact operators \cite[Thm.~2.18]{DP}. Clutching therefore gives only the trivial stabilized bundle on $S^2$, and $\Br_L(S^2)=0$ by \cite[Thm.~2.15(ii)]{DPBrauer}, whereas $\FinBr(S^2)=\Q/\Z$. The finite proof of \eqref{eq:brauer} uses no direct limit over all coprime pairs. Under coordinatewise divisibility those pairs do not form a directed set: $(2,3)$ and $(3,2)$ have no coprime common upper bound. The stabilization of \emph{embedded} algebra bundles with their embeddings retained is a separate construction \cite{Ershov2000,Ershov2008}.

\begin{prop}\label{prop:lowbr}
On connected finite CW complexes of dimension at most three there is a natural isomorphism
\[
                   \FinBr(X)\xrightarrow{\ \theta_X\ } H^2(X;\Q/\Z).
\]
Under this isomorphism the Brauer projection is the Bockstein for
$0\to\Z\to\Q\to\Q/\Z\to0$.
\end{prop}
\begin{proof}
For an algebra $A$ of degree $k$ send $s(A)$ to $H^2(X;\Q/\Z)$ by the coefficient map $[j]\mapsto j/k\bmod\Z$. For the scalar algebra assign zero. Tensoring the special-unitary lifting extensions sends the two central elements $[i]\in\Z/k$ and $[j]\in\Z/m$ to $[mi+kj]\in\Z/(km)$. The resulting classes therefore add in $H^2(X;\Q/\Z)$. Since $s$ vanishes on coprime-embeddable factors, $\theta_X$ is well defined and natural in every dimension. The coefficient map fits into a map of Bockstein sequences whose maps on the two integral terms are the identity on the left and $n\mapsto n/k\in\Q$ in the middle. Thus $\beta\theta_X([A])=\dd(A)$.

In dimension at most three, $c_1:I\to H^2(X;\Z)$ is an isomorphism by the low homotopy groups of $BU$, and $I^2=0$ by the skeletal filtration. The same statements hold rationally. On the kernel in Theorem~\ref{thm:brauer}, a rational unit $[V]/n$, $n=\rk V$, maps under $\theta_X\Phi$ to the coefficient reduction of $c_1(V)/n$. Indeed $s(\End V)=\rho_n c_1(V)$. Thus $\theta_X$ identifies the kernel with $H^2(X;\Q)/\im H^2(X;\Z)$ in the coefficient exact sequence
\[
 0\longrightarrow H^2(X;\Q)/\im H^2(X;\Z)\longrightarrow
 H^2(X;\Q/\Z)\xrightarrow{\ \beta\ } H^3(X;\Z)_{\mathrm{tors}}
 \longrightarrow0 .
\]
The right-hand map between the two exact sequences is the identity, proving the isomorphism. This coefficient sequence already appears in \cite{Ershov2008,v14}; the argument here identifies it with the geometric quotient in these dimensions.
\end{proof}

\begin{prop}\label{prop:moore}
The exact sequence \eqref{eq:brauer} admits no splitting natural under pullback on connected finite CW complexes, even in dimension at most three. On all connected finite CW complexes the functor $\FinBr$ is not representable by a space, and in particular cannot be the degree-zero functor of a generalized cohomology theory.
\end{prop}
\begin{proof}
Write $M_n(q)=S^n\cup_q e^{n+1}$ for the cofiber of a map $S^n\to S^n$ of degree $q\ge2$, and $i_n:S^n\to M_n(q)$ for the inclusion. For $n=2$, positive-degree rational cohomology of $M_2(q)$ vanishes and $H^3(M_2(q);\Z)=\Z/q$. Thus the Bockstein of Proposition~\ref{prop:lowbr} is an isomorphism on $M_2(q)$. The cellular cochains with coefficients in $\Q/\Z$ identify $i_2^*$ on $H^2$ with
\[
         (\Q/\Z)[q]\;\hookrightarrow\;\Q/\Z .
\]
In particular a generator restricts nontrivially to $S^2$, whereas its Brauer class restricts to zero because $H^3(S^2;\Z)=0$. Naturality of a section would force this restriction to be zero.

For $n=4$ the reduced $K$-theory cofiber sequence identifies $\widetilde K^0(M_4(q))$ with the kernel of $q:\Z\to\Z$, since $\widetilde K^{-1}(S^4)=0$. Thus $\widetilde K^0(M_4(q))=0$. Also $H^3(M_4(q);\Z)=0$, so every matrix bundle there is an $\End(V)$ with $[V]=\rk V$. Theorem~\ref{thm:main} now gives $\FinBr(M_4(q))=0$ directly. On $\FinBr(S^4)=\Q/\Z$, the degree-$q$ map induces multiplication by $q$: it multiplies a reduced Bott generator by $q$, and that generator has square zero. A nonzero element of order $q$ thus pulls back to zero along the attaching map. For a functor $[\,\cdot\,,Z]$ any such class on $S^4$ extends over $M_4(q)$ by choosing a null-homotopy on the attaching sphere. The zero class corresponds to a constant map by naturality from a point. An extension of the chosen nonzero class therefore contradicts $\FinBr(M_4(q))=0$. The argument also applies to based representability and gives a direct failure of cofiber exactness for a generalized cohomology theory.
\end{proof}
The statement concerns the quotient of finite matrix bundles by the tensor equivalence of this section. Theories that retain embedding or descent data have a separate representability question.

\section{Groupoid bundles and the family index}\label{sec:fred}
\subsection*{Partial isomorphisms and descent}
Let $G_{k,t}$ have object space $\Gr_{k,t}$ and arrows all $*$-isomorphisms between the corresponding subalgebras. The gauge groupoid of $F_{k,t}\to\Gr_{k,t}$ gives a concrete model:
\[
     G_{k,t}^{(1)}=(F_{k,t}\times F_{k,t})/\PU(k),\qquad
     [g,f]\longmapsto g f^{-1}:\im(f)\longrightarrow\im(g).
\]
The inverse is obtained locally by choosing a frame of the source and transporting it by the given isomorphism. These local inverses agree modulo simultaneous precomposition. Composition and inversion are composition and inversion of algebra maps. Hence $G_{k,t}$ is transitive and its frame bundle implements a Morita equivalence with $\PU(k)$, giving $BG_{k,t}\simeq B\PU(k)$ \cite{Rossi}. Extending the universal $\PU(k)$-bundle along this equivalence gives $EP\times_P F_{k,t}\to BP$ with the groupoid action by postcomposition. The division map between two embeddings $f,g$ of a fixed algebra sends them to the unique arrow $g f^{-1}$. This proves principality and identifies the universal bundle of \cite[Thm.~18]{v14}. A section of the associated groupoid bundle is precisely an embedding, even when the underlying matrix algebra bundle is nontrivial.

If one records the map between commutants as well, an arrow is a pair of isomorphisms of the two tensor factors. There is a unique resulting automorphism of $M_{kt}$. This enlarged groupoid is the action groupoid of $\PU(kt)$ on $\Gr_{k,t}$ and has classifying space $B(\PU(k)\times\PU(t))$. The kernel of the forgetful functor over an object is the automorphism group of that object's commutant. Thus a family of partial isomorphisms extends to ambient automorphisms exactly when isomorphisms of the two commutant bundles can be chosen. For $k,t\ge2$, the identity family of frames on $F_{k,t}$ gives a concrete failure of extension from the constant standard frame: such an extension would be a section of $\PU(kt)\to F_{k,t}$. It would yield a retraction onto the fiber $\PU(t)$ after choosing the section through the identity. On $\pi_3$ the fiber inclusion is multiplication by $k$ on $\Z$ and has no left inverse.

A groupoid cocycle on a cover has maps from the open sets to $\Gr_{k,t}$ and arrows between their restrictions on overlaps, satisfying the usual composition law on triple overlaps. Gluing the image algebras by those arrows gives a global $M_k$-bundle. Conversely, local frames of a global bundle and local embeddings give such a cocycle, up to changes of local frames and refinement. To glue chosen ambient algebras together with the embeddings also requires compatible commutant arrows. This is the local meaning of the two groupoids above.

\subsection*{The multiplicity carried by a frame}
Let $f:X\to F_{k,t}$ and set $(E_f)_x=f_x(e_{11})\C^{kt}$. The isomorphism of Lemma~\ref{lem:mult} in this case is explicit:
\begin{equation}\label{eq:frames}
 \C^k\otimes E_f\longrightarrow\triv{kt},\qquad
 v_i\otimes w\longmapsto f_x(e_{i1})w ,
\end{equation}
where $(v_i)$ is the standard orthonormal basis. Thus a frame family is a rank-$t$ bundle together with a specified unitary trivialization of its $k$-fold sum, modulo isomorphism of such data. Its commutant is $\End(E_f)$. The homogeneous-space fibration $F_{k,t}\to B\PU(t)\to B\PU(kt)$ likewise identifies lifts of a commutant classifying map with null-homotopies of its amplification. Formula~\eqref{eq:frames} identifies the vector bundle that enters the index formula directly from these data.

\begin{prop}\label{prop:fred}
Let $T$ be a norm-continuous Fredholm family on $\C^k\otimes\mathcal H$, with $\mathcal H$ a separable infinite-dimensional Hilbert space. Applying $f_x\otimes\id_{B(\mathcal H)}$ to each $T_x$ gives a Fredholm family on $\C^{kt}\otimes\mathcal H$, with index
\[
                          [E_f]\,\Ind(T)\in K^0(X).
\]
For $t=1$ this gives the line-bundle action; for general $t$, the normalized multiplier $[E_f]/t$ is a unit in $K^0(X)[1/t]$. These normalized multipliers multiply under tensor products of frame families, after passage to a common coefficient ring.
\end{prop}
\begin{proof}
Using \eqref{eq:frames} and permuting the factors, the resulting operator is $T_x\otimes1_{(E_f)_x}$. In the vector-bundle description of the family index, a finite-dimensional stabilization making $T$ surjective can be tensored with $E_f$. Both its kernel class and the stabilizing bundle are thereby multiplied by $[E_f]$. This proves the formula, the same index argument as for the Picard action \cite[Prop.~2.1]{AS}. The localized class $[E_f]/t$ has rank one and nilpotent reduced part. For two frame families the corner projection of their tensor product is the tensor product of their corner projections, which proves the multiplicativity.
\end{proof}

For $t=1$ the family lands in $\PU(k)$ and its line bundle $L$ comes with $L^{\oplus k}\cong\triv{k}$. Tensor order dividing $k$ alone need not supply this trivialization. For example, the complexified tautological line on $\mathbb{RP}^4$ has tensor square one and $c_2(L\oplus L)=c_1(L)^2\ne0$. For $t>1$ the unnormalized operation multiplies numerical ranks by $t$. The normalization of this operation and the action proposed in \cite[§4.3]{v14} are treated below.

\newpage
\subsection*{Stabilized frame actions and higher twists}

\begin{prop}[Integral clutching from a finite frame]\label{prop:frame-unit}
Let $Y$ be a connected finite CW complex, $(k,t)=1$, and $f:Y\to F_{k,t}$. Set $I_Y=\widetilde K^0(Y)$, $I_Y[k]=\{z\in I_Y:kz=0\}$, and $x_f=[E_f]-t$. There is a unique $u_f\in 1+I_Y[k]$ with $t u_f=[E_f]$. It is a unit of integral $K^0(Y)$ given by
\begin{equation}\label{eq:frame-unit}
                u_f=1+a x_f,\qquad at\equiv1\pmod k .
\end{equation}
It determines a $KU$-line on $\Sigma Y$, natural in $Y$. If $\dim Y\le 2t+1$, the units obtained in this way are exactly $1+I_Y[k]$.
\end{prop}
\begin{proof}
Equation~\eqref{eq:frames} gives $E_f^{\oplus k}\cong\triv{kt}$, hence $kx_f=0$. Multiplication by $t$ on $I_Y[k]$ has inverse multiplication by $a$. This proves the formula and uniqueness, even when $K^0(Y)$ also has $t$-torsion. Nilpotence of $I_Y$ makes $u_f$ a unit. Its image after inverting $t$ is the multiplier $[E_f]/t$ of Proposition~\ref{prop:fred}.

Let $KU$ denote the periodic complex $K$-theory spectrum. Its space of units $GL_1(KU)$ is the union of the components of rank $\pm1$ in $\Omega^\infty KU$, with multiplication induced by tensor product. Thus $u_f$ gives a homotopy class of maps $Y\to GL_1(KU)$ of rank one. Clutch two trivial rank-one $KU$-modules over the cones of $\Sigma Y$ by this map, with transition from the negative cone to the positive cone. This gives a $KU$-line in the classification by $BGL_1(KU)$ \cite{ABGHR}. Reversing the convention gives the inverse twist. Pullback preserves the corner and the clutching class.

Conversely, for $z\in I_Y[k]$ choose a rank-$t$ bundle $E$ with $[E]=t(1+z)$, using $\dim Y\le 2t+1$. Then $[E^{\oplus k}]=kt$. Since $k\ge2$ and $\dim Y\le2kt$, cancellation gives a unitary trivialization of $E^{\oplus k}$. By Lemma~\ref{lem:mult} applied to $V=\triv{k}$ this gives a frame with corner $E$. Uniqueness in \eqref{eq:frame-unit} gives $u_f=1+z$.
\end{proof}

\begin{example}[A higher twist of exact order $q$]\label{ex:integral-frame-moore}
Let $q\ge2$, take $\ell\ge2$ coprime to $q$, and put $Y=S^3\cup_q e^4$, $X=\Sigma Y$. The cofiber sequence and Bott periodicity give $\widetilde K^0(Y)=\Z/q$. Products of reduced classes vanish since $Y$ is a suspension. Choose a generator $g$. Proposition~\ref{prop:frame-unit} with $k=q$, $t=\ell$ gives a frame with corner class $\ell+\ell g$ and integral clutching unit $1+g$, since $\dim Y=4\le 2\ell+1$. The normalized clutching group on a suspension of a connected space is the group of rank-one units: changing cone trivializations multiplies a transition class by constant units $\pm1$, and the multiplication of units is homotopy commutative. Thus $(1+g)^j=1+jg$ proves that the twist has exact order $q$. Both $H^1(X;\Z/2)$ and $H^3(X;\Z)$ vanish, so it is a higher twist. Its localization at $\ell$ has unit $[E_f]/\ell=1+g$ as in the prescribed-size calculation \cite{Prescribed}. Compare the related $q=3$ example in \cite[Example~3.5]{DPBrauer}.
\end{example}

Fix coprime $k,\ell\ge2$ and stabilize both matrix sizes. Write
\[
   F_n=F_{k^n,\ell^n},\qquad X_n=\Gr_{k^n,\ell^n},\qquad
   P_n=\PU(k^n)\qquad(n\ge0).
\]
The three spaces $F_0,X_0,BP_0$ are points. The transition on frames is right tensor product with the standard embedding $e:M_k\to M_k\otimes M_\ell$. On images it is tensor product with $\im e$, and on groups it is $u\mapsto u\otimes1_k$. Denote the respective telescopes of $F_n,X_n,BP_n$ by $F_\infty,X_\infty,B_\infty$. The map $p:X_\infty\to B_\infty$ in the homotopy category is induced by the classifying maps of the stagewise frame bundles in Proposition~\ref{prop:grass}.

\begin{theorem}[Action on the tensor-power tower]\label{thm:frame-action}
Tensor product gives $F_\infty$ a grouplike $A_\infty$-space structure up to equivalence. It acts up to coherent homotopy on $X_\infty$, with $p$ as the map to its homotopy quotient. In particular,
\begin{equation}\label{eq:frame-action}
                X_\infty /\!/ F_\infty\ \simeq\ B_\infty
                 =\operatorname{Tel}_n B\PU(k^n).
\end{equation}
Here $/\!/$ denotes the two-sided bar construction of the action.
\end{theorem}
\begin{proof}
Use the Borel model $Q_n=(F_n)_{hP_n}$ for $X_n$. Projection to $BP_n$ has homotopy fiber $F_n$, and $Q_n\to X_n$ is a weak equivalence since precomposition is free and $F_n\to X_n$ is a principal bundle. Use ordered tensor powers and the nerve of the action groupoid to form $Q_n$. The Kronecker product is associative in ordered bases, and tensoring frames intertwines precomposition with the tensor homomorphism $P_m\times P_n\to P_{m+n}$. Thus
\[
\begin{matrix}
 \coprod_n F_n &\longrightarrow& \coprod_n Q_n\\
 \big\downarrow &&\big\downarrow\\
 \mathbb N &\longrightarrow& \coprod_n BP_n
\end{matrix}
\]
is a square of graded topological monoids. The bottom map selects the basepoint in each degree. The square is homotopy cartesian in each degree. Levelwise singular nerves and diagonal realization give well-pointed models throughout.

Each monoid has component monoid $\mathbb N$. On $X_m\times X_n$ the two orders of multiplication differ by an ambient block permutation, joined to the identity through unitaries. On frames use paths to the block permutations in both source and target, acting by inverse precomposition and by postcomposition. On classifying spaces the two tensor homomorphisms differ by an inner conjugation and therefore give homotopic maps. For $\coprod_n Q_n$ use its multiplicative degreewise equivalence to $\coprod_n X_n$. This proves the centrality hypothesis of the group-completion theorem.

For each of the four monoids $M$, the McDuff--Segal theorem identifies the homology of the telescope of its successive components $M_n$ with that of $(\Omega BM)_0$ \cite{MS}. This is the degree-zero component of the usual group-completion comparison: translate the image of degree $n$ by the inverse of the $n$th power of the chosen degree-one point. The four maps preserve the chosen points, so the comparisons form a homotopy-commutative square, functorially up to homotopy.

Here these homology equivalences are weak equivalences. The spaces $BP_n$ and $X_n$ are simply connected, the latter by the special-unitary homogeneous-space description in Section~\ref{sec:frames}. Also $F_n$ is simple. Indeed in its presentation as $U(k^n\ell^n)/(1_{k^n}\otimes U(\ell^n))$ the stabilizer is connected, so each element of $\pi_1 F_n$ comes from a loop in the unitary group. Acting on a representative sphere by that loop gives a homotopy with precisely the required basepoint motion and with identical endpoints. Thus the fundamental group acts trivially on all homotopy groups. Simplicity passes to the telescope, since homotopy groups and this action pass to the sequential limit. The targets $(\Omega BM)_0$ are connected loop spaces and hence simple. The homological Whitehead theorem for simple spaces now applies \cite[Theorem~0.0.2]{MayPonto}.

Taking telescopes in the stagewise homotopy fiber sequences gives another homotopy fiber sequence: compare long exact homotopy sequences using $\pi_j\operatorname{Tel}_n Z_n=\operatorname{colim}_n\pi_j Z_n$ for the based connected towers in question. The lower left telescope is contractible. The group-completion comparisons therefore identify $F_\infty$ with the homotopy fiber of $p$ as a map of grouplike $A_\infty$-spaces.

Deloop $p$ once and put $W=\operatorname{hofib}(BX_\infty\to BB_\infty)$. Both deloopings are simply connected since $X_\infty$ and $B_\infty$ are connected, so $W$ is connected. Looping gives $\Omega W\simeq F_\infty$. The connecting map $B_\infty\to W$ has homotopy fiber $X_\infty$, with fiber inclusion $p$. Pulling back the based path fibration of $W$ gives the action of $\Omega W$ and the asserted homotopy quotient.
\end{proof}

For the enlarged groupoid the kernel of restriction to partial isomorphisms at the standard object is $\PU(\ell^n)$. Stabilization sends $v$ to $v\otimes1_\ell$ and acts by $\ell$ on $\pi_3=\Z$ for $n\ge1$. The telescope of the kernels consequently has $\pi_3=\Z[1/\ell]$. Formula~\eqref{eq:frame-action} identifies the homotopy quotient for the partial-isomorphism groupoids at the level of classifying spaces. It retains the frame data discarded by the tensor equivalence of Section~\ref{sec:quotient}.

On a common parameter space $Y$, frames $f,g$ of degrees $m,n$ have $E_{f\otimes g}\cong E_f\otimes E_g$. Both $u_{f\otimes g}$ and $u_fu_g$ lie in $1+I_Y[k^{m+n}]$ and become $[E_{f\otimes g}]$ after multiplication by $\ell^{m+n}$. Uniqueness in Proposition~\ref{prop:frame-unit} gives $u_{f\otimes g}=u_fu_g$. Since $u_e=1$, stabilization preserves the integral unit. This compatibility is an identity of homotopy classes. Extending the suspension construction to arbitrary principal $F_\infty$-bundles by this prescription requires a coherently multiplicative map $F_\infty\to GL_1(KU)$ realizing these classes, hence a map $BF_\infty\to BGL_1(KU)$. Constructing this map is a separate question from the action in Theorem~\ref{thm:frame-action}.

\subsection*{Acknowledgments}
I thank A. S. Mishchenko and E. V. Troitsky for their support and for many helpful discussions. I would like to express my deep gratitude to Thomas Schick for hospitality and very helpful discussions during my visit to G\"ottingen. I also thank him for pointing out the divisibility of Chern classes on connected covers \cite{Peterson}.
\smallskip

\noindent\emph{Use of generative AI.} OpenAI's ChatGPT assisted in reviewing the earlier preprint, developing and checking arguments, searching the literature, and drafting this revision. The author is responsible for the final text.


\begin{thebibliography}{99}
\bibitem{ABGHR} M. Ando, A. J. Blumberg, D. Gepner, M. J. Hopkins and C. Rezk, \emph{An $\infty$-categorical approach to $R$-line bundles, $R$-module Thom spectra, and twisted $R$-homology}, J. Topol. \textbf{7} (2014), 869--893. \href{https://doi.org/10.1112/jtopol/jtt035}{doi:10.1112/jtopol/jtt035}.
\bibitem{AS} M. Atiyah and G. Segal, \emph{Twisted $K$-theory and cohomology}, \href{https://arxiv.org/abs/math/0510674}{arXiv:math/0510674}.
\bibitem{DP} M. Dadarlat and U. Pennig, \emph{A Dixmier--Douady theory for strongly self-absorbing $C^*$-algebras}, J. Reine Angew. Math. \textbf{718} (2016), 153--181. \href{https://arxiv.org/abs/1302.4468}{arXiv:1302.4468}.
\bibitem{DPBrauer} M. Dadarlat and U. Pennig, \emph{A Dixmier--Douady theory for strongly self-absorbing $C^*$-algebras II: the Brauer group}, J. Noncommut. Geom. \textbf{9} (2015), 1137--1154. \href{https://doi.org/10.4171/JNCG/218}{doi:10.4171/JNCG/218}.
\bibitem{DK} P. Donovan and M. Karoubi, \emph{Graded Brauer groups and $K$-theory with local coefficients}, Publ. Math. IH\'ES \textbf{38} (1970), 5--25. \href{https://doi.org/10.1007/BF02684650}{doi:10.1007/BF02684650}.
\bibitem{Ershov2000} A. V. Ershov, \emph{On the $K$-theory of matrix algebra bundles}, Russian Math. Surveys \textbf{55} (2000), 336--337. \href{https://doi.org/10.1070/RM2000v055n02ABEH000276}{doi:10.1070/RM2000v055n02ABEH000276}.
\bibitem{Ershov2008} A. V. Ershov, \emph{A generalization of the topological Brauer group}, J. K-Theory \textbf{2} (2008), 407--444. \href{https://doi.org/10.1017/is008001006jkt027}{doi:10.1017/is008001006jkt027}.
\bibitem{v14} A. V. Ershov, \emph{Topological obstructions to embedding of a matrix algebra bundle into a trivial one}, 29 August 2010. \href{https://arxiv.org/abs/0807.3544v14}{arXiv:0807.3544v14}.
\bibitem{published} A. V. Ershov, \emph{Obstructions to embeddings of bundles of matrix algebras in a trivial bundle}, Math. Notes \textbf{94} (2013), 482--498. \href{https://doi.org/10.1134/S0001434613090198}{doi:10.1134/S0001434613090198}.
\bibitem{Prescribed} A. V. Ershov, \emph{Prescribed matrix size for UHF Brauer classes on suspensions}, companion manuscript.
\bibitem{HatcherAT} A. Hatcher, \emph{Algebraic Topology}, Cambridge University Press, 2002, Example 2.43 and §3.E.
\bibitem{HatcherVB} A. Hatcher, \emph{Vector Bundles and $K$-Theory}, \href{https://pi.math.cornell.edu/~hatcher/VBKT/VBpage.html}{online monograph}.
\bibitem{Husemoller} D. Husemoller, \emph{Fibre Bundles}, 3rd ed., Graduate Texts in Mathematics 20, Springer, 1994.
\bibitem{MayPonto} J. P. May and K. Ponto, \emph{More Concise Algebraic Topology: Localization, Completion, and Model Categories}, University of Chicago Press, 2012.
\bibitem{MS} D. McDuff and G. Segal, \emph{Homology fibrations and the \textquotedblleft group-completion\textquotedblright\ theorem}, Invent. Math. \textbf{31} (1976), 279--284. \href{https://doi.org/10.1007/BF01403148}{doi:10.1007/BF01403148}.
\bibitem{Peterson} F. P. Peterson, \emph{Some remarks on Chern classes}, Ann. of Math. (2) \textbf{69} (1959), 414--420.
\bibitem{Rossi} C. A. Rossi, \emph{Principal bundles with groupoid structure: local vs. global theory and nonabelian \v{C}ech cohomology}, \href{https://arxiv.org/abs/math/0404449}{arXiv:math/0404449}.
\bibitem{Whitehead} G. W. Whitehead, \emph{Elements of Homotopy Theory}, Graduate Texts in Mathematics 61, Springer, 1978.
\end{thebibliography}
\end{document}